\documentclass{amsart}
\usepackage{amssymb}

\theoremstyle{plain}
\newtheorem{theorem}{Theorem}
\newtheorem{axiom}[theorem]{Axiom}

\newtheorem{conjecture}[theorem]{Conjecture}
\newtheorem{corollary}[theorem]{Corollary}
\newtheorem{exercise}[theorem]{Exercise}
\newtheorem{lemma}[theorem]{Lemma}

\newtheorem{proposition}[theorem]{Proposition}

\theoremstyle{definition}
\newtheorem{definition}[theorem]{Definition}
\newtheorem{remark}[theorem]{Remark}

\newtheorem{example}[theorem]{Example}

\makeatletter
\newcount\@hour\newcount\@minute\chardef\@x10\chardef\@xv60
\def\tcitime{
\def\@time{%
  \@minute\time\@hour\@minute\divide\@hour\@xv
  \ifnum\@hour<\@x 0\fi\the\@hour:%
  \multiply\@hour\@xv\advance\@minute-\@hour
  \ifnum\@minute<\@x 0\fi\the\@minute
  }}%

\@ifundefined{hyperref}{}{}

\@ifundefined{qExtProgCall}{\def\qExtProgCall#1#2#3#4#5#6{\relax}}{}
\def\QCTOpt[#1]#2{%
  \def\QCTOptB{#1}
  \def\QCTOptA{#2}
}
\def\QCTNOpt#1{%
  \def\QCTOptA{#1}
  \let\QCTOptB\empty
}
\def\Qct{%
  \@ifnextchar[{%
    \QCTOpt}{\QCTNOpt}
}
\def\QCBOpt[#1]#2{%
  \def\QCBOptB{#1}
  \def\QCBOptA{#2}
}
\def\QCBNOpt#1{%
  \def\QCBOptA{#1}
  \let\QCBOptB\empty
}
\def\Qcb{%
  \@ifnextchar[{%
    \QCBOpt}{\QCBNOpt}
}
\def\PrepCapArgs{%
  \ifx\QCBOptA\empty
    \ifx\QCTOptA\empty
      {}%
    \else
      \ifx\QCTOptB\empty
        {\QCTOptA}%
      \else
        [\QCTOptB]{\QCTOptA}%
      \fi
    \fi
  \else
    \ifx\QCBOptA\empty
      {}%
    \else
      \ifx\QCBOptB\empty
        {\QCBOptA}%
      \else
        [\QCBOptB]{\QCBOptA}%
      \fi
    \fi
  \fi
}
\newcount\GRAPHICSTYPE
\GRAPHICSTYPE=\z@
\def\GRAPHICSPS#1{%
 \ifcase\GRAPHICSTYPE
   \special{ps: #1}%
 \or
   \special{language "PS", include "#1"}%
 \fi
}%
\def\graffile#1#2#3#4{%
    \leavevmode
    \raise -#4 \BOXTHEFRAME{%
        \hbox to #2{\raise #3\hbox to #2{\null #1\hfil}}}%
}%
\def\draftbox#1#2#3#4{%
 \leavevmode\raise -#4 \hbox{%
  \frame{\rlap{\protect\tiny #1}\hbox to #2%
   {\vrule height#3 width\z@ depth\z@\hfil}%
  }%
 }%
}%
\newcount\draft
\draft=\z@

\newif\ifwasdraft
\wasdraftfalse

\def\GRAPHIC#1#2#3#4#5{%
 \ifnum\draft=\@ne\draftbox{#2}{#3}{#4}{#5}%
  \else\graffile{#1}{#3}{#4}{#5}%
  \fi
 }%
\def\addtoLaTeXparams#1{%
    \edef\LaTeXparams{\LaTeXparams #1}}%
\newif\ifBoxFrame \BoxFramefalse
\newif\ifOverFrame \OverFramefalse
\newif\ifUnderFrame \UnderFramefalse

\def\BOXTHEFRAME#1{%
   \hbox{%
      \ifBoxFrame
         \frame{#1}%
      \else
         {#1}%
      \fi
   }%
}

\def\doFRAMEparams#1{\BoxFramefalse\OverFramefalse\UnderFramefalse\readFRAMEparams#1\end}%
\def\readFRAMEparams#1{%
 \ifx#1\end%
  \let\next=\relax
  \else
  \ifx#1i\dispkind=\z@\fi
  \ifx#1d\dispkind=\@ne\fi
  \ifx#1f\dispkind=\tw@\fi
  \ifx#1t\addtoLaTeXparams{t}\fi
  \ifx#1b\addtoLaTeXparams{b}\fi
  \ifx#1p\addtoLaTeXparams{p}\fi
  \ifx#1h\addtoLaTeXparams{h}\fi
  \ifx#1X\BoxFrametrue\fi
  \ifx#1O\OverFrametrue\fi
  \ifx#1U\UnderFrametrue\fi
  \ifx#1w
    \ifnum\draft=1\wasdrafttrue\else\wasdraftfalse\fi
    \draft=\@ne
  \fi
  \let\next=\readFRAMEparams
  \fi
 \next
 }%
\def\IFRAME#1#2#3#4#5#6{%
      \bgroup
      \let\QCTOptA\empty
      \let\QCTOptB\empty
      \let\QCBOptA\empty
      \let\QCBOptB\empty
      #6%
      \parindent=0pt%
      \leftskip=0pt
      \rightskip=0pt
      \setbox0 = \hbox{\QCBOptA}%
      \@tempdima = #1\relax
      \ifOverFrame
          \typeout{This is not implemented yet}%
          \show\HELP
      \else
         \ifdim\wd0>\@tempdima
            \advance\@tempdima by \@tempdima
            \ifdim\wd0 >\@tempdima
               \textwidth=\@tempdima
               \setbox1 =\vbox{%
                  \noindent\hbox to \@tempdima{\hfill\GRAPHIC{#5}{#4}{#1}{#2}{#3}\hfill}\\%
                  \noindent\hbox to \@tempdima{\parbox[b]{\@tempdima}{\QCBOptA}}%
               }%
               \wd1=\@tempdima
            \else
               \textwidth=\wd0
               \setbox1 =\vbox{%
                 \noindent\hbox to \wd0{\hfill\GRAPHIC{#5}{#4}{#1}{#2}{#3}\hfill}\\%
                 \noindent\hbox{\QCBOptA}%
               }%
               \wd1=\wd0
            \fi
         \else
            \ifdim\wd0>0pt
              \hsize=\@tempdima
              \setbox1 =\vbox{%
                \unskip\GRAPHIC{#5}{#4}{#1}{#2}{0pt}%
                \break
                \unskip\hbox to \@tempdima{\hfill \QCBOptA\hfill}%
              }%
              \wd1=\@tempdima
           \else
              \hsize=\@tempdima
              \setbox1 =\vbox{%
                \unskip\GRAPHIC{#5}{#4}{#1}{#2}{0pt}%
              }%
              \wd1=\@tempdima
           \fi
         \fi
         \@tempdimb=\ht1
         \advance\@tempdimb by \dp1
         \advance\@tempdimb by -#2%
         \advance\@tempdimb by #3%
         \leavevmode
         \raise -\@tempdimb \hbox{\box1}%
      \fi
      \egroup%
}%
\def\DFRAME#1#2#3#4#5{%
 \begin{center}
     \let\QCTOptA\empty
     \let\QCTOptB\empty
     \let\QCBOptA\empty
     \let\QCBOptB\empty
     \ifOverFrame 
        #5\QCTOptA\par
     \fi
     \GRAPHIC{#4}{#3}{#1}{#2}{\z@}
     \ifUnderFrame 
        \nobreak\par #5\QCBOptA
     \fi
 \end{center}%
 }%
\def\FFRAME#1#2#3#4#5#6#7{%
 \begin{figure}[#1]%
  \let\QCTOptA\empty
  \let\QCTOptB\empty
  \let\QCBOptA\empty
  \let\QCBOptB\empty
  \ifOverFrame
    #4
    \ifx\QCTOptA\empty
    \else
      \ifx\QCTOptB\empty
        \caption{\QCTOptA}%
      \else
        \caption[\QCTOptB]{\QCTOptA}%
      \fi
    \fi
    \ifUnderFrame\else
      \label{#5}%
    \fi
  \else
    \UnderFrametrue%
  \fi
  \begin{center}\GRAPHIC{#7}{#6}{#2}{#3}{\z@}\end{center}%
  \ifUnderFrame
    #4
    \ifx\QCBOptA\empty
      \caption{}%
    \else
      \ifx\QCBOptB\empty
        \caption{\QCBOptA}%
      \else
        \caption[\QCBOptB]{\QCBOptA}%
      \fi
    \fi
    \label{#5}%
  \fi
  \end{figure}%
 }%
\newcount\dispkind%

\def\makeactives{
  \catcode`\"=\active
  \catcode`\;=\active
  \catcode`\:=\active
  \catcode`\'=\active
  \catcode`\~=\active
}
\bgroup
   \makeactives
   \gdef\activesoff{%
      \def"{\string"}
      \def;{\string;}
      \def:{\string:}
      \def'{\string'}
      \def~{\string~}
    }
\egroup

\def\FRAME#1#2#3#4#5#6#7#8{%
 \bgroup
 \@ifundefined{bbl@deactivate}{}{\activesoff}
 \ifnum\draft=\@ne
   \wasdrafttrue
 \else
   \wasdraftfalse%
 \fi
 \def\LaTeXparams{}%
 \dispkind=\z@
 \def\LaTeXparams{}%
 \doFRAMEparams{#1}%
 \ifnum\dispkind=\z@\IFRAME{#2}{#3}{#4}{#7}{#8}{#5}\else
  \ifnum\dispkind=\@ne\DFRAME{#2}{#3}{#7}{#8}{#5}\else
   \ifnum\dispkind=\tw@
    \edef\@tempa{\noexpand\FFRAME{\LaTeXparams}}%
    \@tempa{#2}{#3}{#5}{#6}{#7}{#8}%
    \fi
   \fi
  \fi
  \ifwasdraft\draft=1\else\draft=0\fi{}%
  \egroup
 }%
\def\TEXUX#1{"texux"}

\def\limfunc#1{\mathop{\rm #1}}%

\long\def\QQQ#1#2{%
     \long\expandafter\def\csname#1\endcsname{#2}}%
\@ifundefined{QTP}{\def\QTP#1{}}{}
\@ifundefined{QEXCLUDE}{\def\QEXCLUDE#1{}}{}
\@ifundefined{Qlb}{}{}
\@ifundefined{Qlt}{}{}
\long\def\QQA#1#2{}%
\def\QTR#1#2{{\csname#1\endcsname #2}}
\def\EXPAND#1[#2]#3{}%
\def\NOEXPAND#1[#2]#3{}%
\def\LaTeXparent#1{}%
\def\ChildStyles#1{}%
\def\ChildDefaults#1{}%
\def\QTagDef#1#2#3{}%
\@ifundefined{StyleEditBeginDoc}{}{}
\def\QQfnmark#1{\footnotemark}

\@ifundefined{INDEX}{\def\INDEX#1#2{}{}}{}%
\@ifundefined{SUBINDEX}{\def\SUBINDEX#1#2#3{}{}{}}{}%
\@ifundefined{initial}%
   {\def\initial#1{\bigbreak{\raggedright\large\bf #1}\kern 2\p@\penalty3000}}%
   {}%
\@ifundefined{entry}{}{}%
\@ifundefined{primary}{}{}%
\@ifundefined{secondary}{}{}%
\@ifundefined{ZZZ}{}{\makeatletter\input gnuindex.sty\makeatother\makeindex\makeatletter}%
\@ifundefined{abstract}{%
 \def\abstract{%
  \if@twocolumn
   \section*{Abstract (Not appropriate in this style!)}%
   \else \small 
   \begin{center}{\bf Abstract\vspace{-.5em}\vspace{\z@}}\end{center}%
   \quotation 
   \fi
  }%
 }{%
 }%
\@ifundefined{endabstract}{\def\endabstract
  {\if@twocolumn\else\endquotation\fi}}{}%
\@ifundefined{maketitle}{\def\maketitle#1{}}{}%
\@ifundefined{affiliation}{\def\affiliation#1{}}{}%
\@ifundefined{proof}{}{}%
\@ifundefined{endproof}{}{}%
\@ifundefined{newfield}{\def\newfield#1#2{}}{}%
\@ifundefined{chapter}{\def\chapter#1{\par(Chapter head:)#1\par }%
 \newcount\c@chapter}{}%
\@ifundefined{part}{\def\part#1{\par(Part head:)#1\par }}{}%
\@ifundefined{section}{\def\section#1{\par(Section head:)#1\par }}{}%
\@ifundefined{subsection}{\def\subsection#1%
 {\par(Subsection head:)#1\par }}{}%
\@ifundefined{subsubsection}{\def\subsubsection#1%
 {\par(Subsubsection head:)#1\par }}{}%
\@ifundefined{paragraph}{\def\paragraph#1%
 {\par(Subsubsubsection head:)#1\par }}{}%
\@ifundefined{subparagraph}{\def\subparagraph#1%
 {\par(Subsubsubsubsection head:)#1\par }}{}%
\@ifundefined{therefore}{}{}%
\@ifundefined{backepsilon}{}{}%
\@ifundefined{yen}{}{}%
\@ifundefined{registered}{%
   \def\registered{\relax\ifmmode{}\r@gistered
                    \else$\m@th\r@gistered$\fi}%
 \def\r@gistered{^{\ooalign
  {\hfil\raise.07ex\hbox{$\scriptstyle\rm\text{R}$}\hfil\crcr
  \mathhexbox20D}}}}{}%
\@ifundefined{Eth}{}{}%
\@ifundefined{eth}{}{}%
\@ifundefined{Thorn}{}{}%
\@ifundefined{thorn}{}{}%
\@ifundefined{degree}{}{}%
\newdimen\theight
\def\Column{%
 \vadjust{\setbox\z@=\hbox{\scriptsize\quad\quad tcol}%
  \theight=\ht\z@\advance\theight by \dp\z@\advance\theight by \lineskip
  \kern -\theight \vbox to \theight{%
   \rightline{\rlap{\box\z@}}%
   \vss
   }%
  }%
 }%
\def\qed{%
 \ifhmode\unskip\nobreak\fi\ifmmode\ifinner\else\hskip5\p@\fi\fi
 \hbox{\hskip5\p@\vrule width4\p@ height6\p@ depth1.5\p@\hskip\p@}%
 }%
\def\miss{\hbox{\vrule height2\p@ width 2\p@ depth\z@}}%
\def\tcol#1{{\baselineskip=6\p@ \vcenter{#1}} \Column}  %
\def\newfmtname{LaTeX2e}
\def\chkcompat{%
   \if@compatibility
   \else
     \usepackage{latexsym}
   \fi
}

\ifx\fmtname\newfmtname
  \DeclareOldFontCommand{\rm}{\normalfont\rmfamily}{\mathrm}
  \DeclareOldFontCommand{\sf}{\normalfont\sffamily}{\mathsf}
  \DeclareOldFontCommand{\tt}{\normalfont\ttfamily}{\mathtt}
  \DeclareOldFontCommand{\bf}{\normalfont\bfseries}{\mathbf}
  \DeclareOldFontCommand{\it}{\normalfont\itshape}{\mathit}
  \DeclareOldFontCommand{\sl}{\normalfont\slshape}{\@nomath\sl}
  \DeclareOldFontCommand{\sc}{\normalfont\scshape}{\@nomath\sc}
  \chkcompat
\fi

\def\alpha{{\Greekmath 010B}}%
\def\beta{{\Greekmath 010C}}%
\def\gamma{{\Greekmath 010D}}%
\def\delta{{\Greekmath 010E}}%
\def\epsilon{{\Greekmath 010F}}%
\def\zeta{{\Greekmath 0110}}%
\def\eta{{\Greekmath 0111}}%
\def\theta{{\Greekmath 0112}}%
\def\iota{{\Greekmath 0113}}%
\def\kappa{{\Greekmath 0114}}%
\def\lambda{{\Greekmath 0115}}%
\def\mu{{\Greekmath 0116}}%
\def\nu{{\Greekmath 0117}}%
\def\xi{{\Greekmath 0118}}%
\def\pi{{\Greekmath 0119}}%
\def\rho{{\Greekmath 011A}}%
\def\sigma{{\Greekmath 011B}}%
\def\tau{{\Greekmath 011C}}%
\def\upsilon{{\Greekmath 011D}}%
\def\phi{{\Greekmath 011E}}%
\def\chi{{\Greekmath 011F}}%
\def\psi{{\Greekmath 0120}}%
\def\omega{{\Greekmath 0121}}%
\def\varepsilon{{\Greekmath 0122}}%
\def\vartheta{{\Greekmath 0123}}%
\def\varpi{{\Greekmath 0124}}%
\def\varrho{{\Greekmath 0125}}%
\def\varsigma{{\Greekmath 0126}}%
\def\varphi{{\Greekmath 0127}}%

\def\nabla{{\Greekmath 0272}}
\def\FindBoldGroup{%
   {\setbox0=\hbox{$\mathbf{x\global\edef\theboldgroup{\the\mathgroup}}$}}%
}

\def\Greekmath#1#2#3#4{%
    \if@compatibility
        \ifnum\mathgroup=\symbold
           \mathchoice{\mbox{\boldmath$\displaystyle\mathchar"#1#2#3#4$}}%
                      {\mbox{\boldmath$\textstyle\mathchar"#1#2#3#4$}}%
                      {\mbox{\boldmath$\scriptstyle\mathchar"#1#2#3#4$}}%
                      {\mbox{\boldmath$\scriptscriptstyle\mathchar"#1#2#3#4$}}%
        \else
           \mathchar"#1#2#3#4%
        \fi 
    \else 
        \FindBoldGroup
        \ifnum\mathgroup=\theboldgroup 
           \mathchoice{\mbox{\boldmath$\displaystyle\mathchar"#1#2#3#4$}}%
                      {\mbox{\boldmath$\textstyle\mathchar"#1#2#3#4$}}%
                      {\mbox{\boldmath$\scriptstyle\mathchar"#1#2#3#4$}}%
                      {\mbox{\boldmath$\scriptscriptstyle\mathchar"#1#2#3#4$}}%
        \else
           \mathchar"#1#2#3#4%
        \fi     	    
	  \fi}

\newif\ifGreekBold  \GreekBoldfalse
\let\SAVEPBF=\pbf
\def\pbf{\GreekBoldtrue\SAVEPBF}%

\@ifundefined{theorem}{\newtheorem{theorem}{Theorem}}{}
\@ifundefined{lemma}{\newtheorem{lemma}[theorem]{Lemma}}{}
\@ifundefined{corollary}{\newtheorem{corollary}[theorem]{Corollary}}{}
\@ifundefined{conjecture}{}{}
\@ifundefined{proposition}{\newtheorem{proposition}[theorem]{Proposition}}{}
\@ifundefined{axiom}{}{}
\@ifundefined{remark}{}{}
\@ifundefined{example}{}{}
\@ifundefined{exercise}{}{}
\@ifundefined{definition}{}{}

\@ifundefined{mathletters}{%
  \newcounter{equationnumber}  
  \def\mathletters{%
     \addtocounter{equation}{1}
     \edef\@currentlabel{\theequation}%
     \setcounter{equationnumber}{\c@equation}
     \setcounter{equation}{0}%
     \edef\theequation{\@currentlabel\noexpand\alph{equation}}%
  }
  
}{}

\@ifundefined{BibTeX}{%
    \def\BibTeX{{\rm B\kern-.05em{\sc i\kern-.025em b}\kern-.08em
                 T\kern-.1667em\lower.7ex\hbox{E}\kern-.125emX}}}{}%
\@ifundefined{AmS}%
    {\def\AmS{{\protect\usefont{OMS}{cmsy}{m}{n}%
                A\kern-.1667em\lower.5ex\hbox{M}\kern-.125emS}}}{}%
\@ifundefined{AmSTeX}{}{}%
\begin{document}
\title{Whitehead modules over large principal ideal domains}
\author{Paul C. Eklof}
\thanks{First author partially supported by NSF DMS 98-03126.}
\address[Eklof]{Math Dept, UCI\\
Irvine, CA 92697-3875}
\author{Saharon Shelah}
\thanks{Second author supported by the German-Israeli Foundation for Scientific
Research \& Development. Publication 752.}
\address[Shelah]{Institute of Mathematics, Hebrew University\\
Jerusalem 91904, Israel}
\date{\today}

\begin{abstract}
We consider the Whitehead problem for principal ideal domains of large size.
It is proved, in ZFC, that some p.i.d.'s of size $\geq \aleph _{2}$ have
non-free Whitehead modules even though they are not complete discrete
valuation rings.
\end{abstract}

\maketitle

A module $M$ is a \textit{Whitehead module} if $\limfunc{Ext}_{R}^{1}(M,R)=0$%
. The second author proved that the problem of whether every Whitehead $\mathbb{%
Z}$-module is free is independent of ZFC + GCH (cf. \cite{Sh74}, \cite{Sh75}%
, \cite{Sh80}). This was extended in \cite{BFS} to modules over principal
ideal domains of cardinality at most $\aleph _{1}$. Here we consider the
Whitehead problem for modules over principal ideal domains (p.i.d.'s) of
cardinality $>\aleph _{1}$.

If $R$ is any p.i.d.\ which is not a complete discrete valuation ring, then
an $R$-module of countable rank is Whitehead if and only if it is free (cf. 
\cite{GKW}). On the other hand, if $R$ is a complete discrete valuation
ring, then it is cotorsion and hence every torsion-free $R$-module is a
Whitehead module (cf. \cite[XII.1.17]{EM}).

It will be convenient to decree that a field is not a p.i.d.\ and to use the
term ``slender'' to designate a p.i.d.\ which is not a complete discrete
valuation ring, or equivalently, is not cotorsion (cf. \cite[III.2.9]{EM}).
We will say that a module is $\kappa $-generated if it is generated by a
subset of size $\leq $ $\kappa $ and that it is $\kappa $-free if every
submodule generated by $<\kappa $ elements is free. (Note that, by
Pontryagin's Criterion and induction on $\kappa $, every $\aleph _{1}$-free
module which has rank $\leq \kappa $ is $\kappa $-generated.)

An argument due to the second author (cf. \cite{Sh80} or \cite{T}) shows
that it is consistent with ZFC + GCH that for any p.i.d.\ $R$ (of arbitrary
size), there are Whitehead $R$-modules of rank $\geq |R|$ which are not free.

If the p.i.d.\ $R$ is slender and has cardinality at most $\aleph _{1}$, the
Axiom of Constructibility (V = L) implies that every Whitehead $R$-module is
free (cf. \cite{BFS}). Our main result is that the story is different for
p.i.d.'s of larger size. We will prove the following theorems in ZFC.

\begin{theorem}
\label{main1}There is a slender p.i.d.\ $R$ of cardinality $2^{\aleph _{1}}$
such that every $\aleph _{1}$-free $\aleph _{1}$-generated $R$-module is a
Whitehead module. Hence there are non-free Whitehead $R$-modules which are $%
\aleph _{1}$-generated.
\end{theorem}

\begin{theorem}
\label{main2}There is a p.i.d.\ $R$ of cardinality $\aleph _{2}$ such that
an $\aleph _{1}$-generated $R$-module is Whitehead only if it is free.
\end{theorem}

Assuming V = L and using the existing theory (cf. \cite{BFS}) one easily
obtains the following:

\begin{corollary}
\label{cor}(V = L) There are principal ideal domains\ $R_{1}$ and $R_{2}$
each of cardinality $\aleph _{2}$ and non-slender such that:

(1) an $R_{1}$-module $M$ (of arbitrary cardinality) is Whitehead if and
only if $M$ is the union of a continuous chain, $M=\bigcup_{\alpha <\lambda
}M_{\alpha }$ for some $\lambda $, such that for all $\alpha <\lambda $, $%
M_{\alpha +1}/M_{\alpha }$ is $\aleph _{1}$-free and $\aleph _{1}$-generated;

(2) an $R_{2}$-module $M$ (of arbitrary cardinality) is Whitehead if and
only if $M$ is free. \qed  
\end{corollary}

The theorems can be generalized to other cardinals: see Theorems \ref{1b}
and \ref{2b} at the end of the sections.

\section{Proof of Theorem \ref{main1}}

The ring $R$ in Theorem \ref{main1} will be constructed by a transfinite
induction so that for every module $F/K$ ($F$ free) which is $\aleph _{1}$%
-free and $\aleph _{1}$-generated, $\limfunc{Ext}(F/K,R)=0$, i.e., every
homomorphism from $K$ to $R$ extends to a homomorphism from $F$ to $R$. The
following proposition provides the inductive step.

\begin{proposition}
\label{rngext}Let $R$ be a local slender p.i.d. with maximal ideal $pR$, \
and let $K\subseteq F$ be free $R$-modules of rank $\aleph _{1}$ such that $%
F/K$ is $\aleph _{1}$-free. Let $\psi :K\rightarrow R$ be an $R$%
-homomorphism. Then there is a local slender p.i.d.\ $R^{+}$ containing $R$
as subring, with maximal ideal $pR^{+}$ and of cardinality $=|R|+\aleph _{1}$
such that the $R^{+}$-homomorphism $1_{R^{+}}\otimes _{R}$ $\psi
:R^{+}\otimes _{R}K\rightarrow R^{+}\otimes _{R}R$ extends to an $R^{+}$%
-homomorphism $\varphi :R^{+}\otimes _{R}F\rightarrow R^{+}\otimes _{R}R$.
\end{proposition}

\noindent \textsc{Proof}. Write $F=\bigcup_{\alpha <\omega _{1}}F_{\alpha }$
as a continuous union of submodules of countable rank with $F_{0}=0$. For
each $\alpha <\omega _{1}$, $F_{\alpha }+K/K$ is free; let $\{b_{i}^{\alpha
}:i\in I_{\alpha }\}$ be a linearly independent subset of $F_{\alpha }$ such
that $\{b_{i}^{\alpha }+K:i\in I_{\alpha }\}$ is a basis of $F_{\alpha }+K/K$%
. ($I_{0}=\emptyset $.) Then for all $\alpha <\beta <\omega _{1}$ and all $%
i\in I_{\alpha }$, $b_{i}^{\alpha }=\sum_{j\in I_{\beta }}r_{i,j}^{\alpha
,\beta }b_{j}^{\beta }+k_{i}^{\alpha ,\beta }$ for some unique $%
r_{i,j}^{\alpha ,\beta }\in R$ (which are equal to $0$ for almost all $j$)
and $k_{i}^{\alpha ,\beta }\in K$. Let $s_{i}^{\alpha ,\beta }=\psi
(k_{i}^{\alpha ,\beta })$.

We claim that there is a local slender p.i.d.\ $R^{+}$ of cardinality $%
=|R|+\aleph _{1}$ containing $R$ as subring and with maximal ideal $pR^{+}$
and elements $x_{i}^{\alpha }$ $\in R^{+}$ ($\alpha <\omega _{1}$, $i\in
I_{\alpha }$) such that $x_{i}^{\alpha }=\sum_{j\in I_{\beta
}}r_{i,j}^{\alpha ,\beta }x_{j}^{\beta }+s_{i}^{\alpha ,\beta }$ for all $%
\alpha <\beta <\omega _{1}$ and $i\in I_{\alpha }$. Supposing this for the
moment, let us finish the proof. Clearly $\{b_{i}^{\alpha }:\alpha <\omega
_{1}$, $i\in I_{\alpha }\}\cup K$ generates $R^{+}\otimes _{R}F$ as $R^{+}$%
-module. Define $\varphi $ extending $1_{R^{+}}\otimes _{R}$ $\psi $ by $%
\varphi (1\otimes b_{i}^{\alpha })=x_{i}^{\alpha }\otimes 1$. We must check
that this is well-defined. For this it suffices to prove that $\varphi
(1\otimes b_{i}^{\alpha })=\sum_{j\in I_{\beta }}r_{i,j}^{\alpha ,\beta
}\varphi (1\otimes b_{j}^{\beta })+(1\otimes \psi )(1\otimes k_{i}^{\alpha
,\beta })$ for all $\alpha <\beta <\omega _{1}$ and $i\in I_{\alpha }$. But
this is implied by the assumption that $x_{i}^{\alpha }=\sum_{j\in I_{\beta
}}r_{i,j}^{\alpha ,\beta }x_{j}^{\beta }+s_{i}^{\alpha ,\beta }$. 

So it remains to define $R^{+}$. Let $R^{0}=R$ and for $0<$ $\alpha <\omega
_{1}$, let $R^{\alpha }=R[\{x_{i}^{\alpha }:i\in I_{\alpha }\}]$, the
polynomial ring over $R$ in the commuting indeterminates $x_{i}^{\alpha }$, $%
i\in I_{\alpha }$. For $\alpha <\beta <\omega _{1}$, let $\pi _{\beta
}^{\alpha }:R^{\alpha }\rightarrow R^{\beta }$ be the ring homomorphism
which is the identity on $R$ and takes $x_{i}^{\alpha }$ to $\sum_{j\in
I_{\beta }}r_{i,j}^{\alpha ,\beta }x_{j}^{\beta }+s_{i}^{\alpha ,\beta }$.
It is easy to check, using the fact that the $\{b_{i}^{\gamma }:i\in
I_{\gamma }\}$ are linearly independent, that $\pi _{\gamma }^{\beta }\circ
\pi _{\beta }^{\alpha }=\pi _{\gamma }^{\alpha }$ whenever $\alpha <\beta
<\gamma <\omega _{1}$. Let $R^{\prime }$ with maps $\pi ^{\alpha }:R^{\alpha
}\rightarrow R^{\prime }$ be the direct limit of this $\aleph _{1}$-directed
system of homomorphisms. Clearly each $R^{\alpha }$ is a unique
factorization domain such that $p$ is prime in $R^{\alpha }$. Since the
system is directed, $R^{\prime }$ is an integral domain and $p$ is prime in $%
R^{\prime }$. Moreover, since the system is $\aleph _{1}$-directed, $%
\bigcap_{n\in \omega }p^{n}R^{\prime }=0$ since the same is true in each $%
R^{\alpha }$. If $\{a_{n}:n\in \omega \}$ is a Cauchy sequence in $R$ which
does not have a limit (in the $p$-adic topology), then $\{t\pi _{\alpha
}^{0}(a_{n}):n\in \omega \}$ does not have a limit in the $p$-adic topology
on $R^{\alpha }$ for all $t\in R^{\alpha }-pR^{\alpha }$. Hence, by the $%
\aleph _{1}$-directedness, the same holds for $\{\pi ^{0}(a_{n}):n\in \omega
\}$ in $R^{\prime }$.

Finally, let $R^{+}$ be the localization of $R^{\prime }$ at the prime $p$.
We appeal to the following elementary Lemma to finish. \qed  

\medskip

\begin{lemma}
\label{local}Suppose $R^{\prime }$ is an integral domain with a prime $p$
such that $\bigcap_{n\in \omega }p^{n}R^{\prime }=0$. Then the localization $%
R_{(p)}^{\prime }$ of $R^{\prime }$ at $p$ is a p.i.d.
\end{lemma}

\noindent \textsc{Proof}. Given a non-zero proper ideal $I$ of $%
R_{(p)}^{\prime }$, let $I^{\prime }=I\cap R^{\prime }$ ($=\{r\in R^{\prime
}:\frac{r}{1}\in I\}$). Let $m$ be minimal such that $I^{\prime }\cap
(p^{m}R^{\prime }-p^{m+1}R^{\prime })\neq \emptyset $. Clearly $m$ exists,
by hypothesis and since $I^{\prime }$ is non-zero. We claim that $%
I=p^{m}R_{(p)}^{\prime }$. Let $a\in I^{\prime }\cap (p^{m}R^{\prime
}-p^{m+1}R^{\prime })$; then $a=p^{m}r$ for some $r\in R^{\prime }$ and $%
r\notin pR^{\prime }$; so $r$ is a unit in $R_{(p)}^{\prime }$ and thus $%
p^{m}\in I$. Now for any non-zero $\frac{b}{t}\in I$, $b\in I^{\prime
}-\{0\} $ so $b\in I^{\prime }\cap (p^{n}R^{\prime }-p^{n+1}R^{\prime })$
for some $n\geq m$. Thus $b=p^{n}c$ for some $c\in R^{\prime }$ and $n\geq m$%
. But then $\frac{b}{t}=p^{n}\frac{c}{t}\in p^{m}R_{(p)}^{\prime }$.
Therefore $I=p^{m}R_{(p)}^{\prime }$. \qed  

\medskip

\noindent \textbf{\ \noindent Proof of Theorem \ref{main1}. Let }$\lambda
=2^{\aleph _{1}}$. We  define a ring $R$ on the set $\lambda $ which is the
union of a continuous chain of rings $R_{\nu }$ ($\nu <\lambda $) such that
for each $\nu <\lambda $, $R_{\nu +1}$ is of the form $(R_{\nu })^{+}$ for
some quadruple $(R_{\nu },K_{\nu },F_{\nu },\psi _{\nu })$ satisfying the
hypotheses of the Proposition. We begin, for example, with $R_{0}=\mathbb{Z}%
_{(p)}$. It is easy to see that $R$ is a local p.i.d.\ with prime $p$.
Moreover, the proof of the Proposition shows that a witnessing Cauchy
sequence to the incompleteness of $R_{0}$ is preserved at each stage and
therefore also in $R$ since $\omega _{1}$ has cofinality $>\omega $. Because 
$\lambda ^{\aleph _{1}}=\lambda $, we can choose the enumeration of
quadruples $(R_{\nu },K_{\nu },F_{\nu },\psi _{\nu })$ such that for every $%
\aleph _{1}$-generated $\aleph _{1}$-free $R$-module $F/K$ (where $%
K\subseteq F$ are free $R$-modules) and every $R$-homomorphism $\psi
:K\rightarrow R$, there is a $\nu <\lambda $ such that $R\otimes _{R_{\nu
}}F_{\nu }$ is isomorphic to $F$ under an isomorphism which takes $R\otimes
_{R_{\nu }}K_{\nu }$ to $K$ and identifies $1_{R}\otimes _{R_{\nu }}\psi
_{\nu }$ with $\psi $ under the natural isomorphism of $R\otimes _{R_{\nu
}}R_{\nu }$ with $R$. (Note that $K\subseteq F$ and $\psi $ can each be
completely described by a sequence of $\aleph _{1}$ elements of $R=\lambda $%
.) \qed  

\medskip

By using a direct system indexed by the countable rank submodules of $F/K$
in the proof of the Proposition, we can prove the following more general
version of the theorem. Part (1) of Corollary \ref{cor} can be
correspondingly generalized.

\begin{theorem}
\label{1b}For any cardinal $\kappa \geq \aleph _{1}$, there is a local
slender p.i.d.\ $R$ of cardinality $2^{\kappa }$ such that every $\aleph _{1}
$-free $\kappa $-generated $R$-module is a Whitehead module. \qed  
\end{theorem}

\section{Proof of Theorem \ref{main2}}

Let $R$ be the polynomial ring $\mathcal{F}[X]$ where $\mathcal{F}=\mathbb{Q}%
(\{t_{\nu }:\nu <\omega _{2}\})$ and $\{t_{\nu }:\nu <\omega _{2}\}$ is an
algebraically independent set.

Let $A$ be an $\aleph _{1}$-generated $\aleph _{1}$-free $R$-module which is
not free and let $A=\bigcup_{\alpha <\omega _{1}}A_{\alpha }$ be an $\aleph
_{1}$-filtration of $A$. Then there is a stationary set $S$ of limit
ordinals such that for $\gamma \in S$, $A_{\gamma +1}/A_{\gamma }$ is not
free. Without loss of generality we can assume that there is a $d\in \omega $
such that for all $\gamma \in S$, $A_{\gamma +1}/A_{\gamma }$ is of rank $%
d+1 $ and not free but every submodule of rank $\leq d$ is free$.$ (Note
that we allow $A_{\alpha +1}/A_{\alpha }$ to be non-free for $\alpha \notin
S $.) Thus $A_{\gamma +1}/A_{\gamma }$ is isomorphic to $F_{\gamma }^{\prime
}/K_{\gamma }^{\prime }$ where $F_{\gamma }^{\prime }$ is free on $%
\{y_{\gamma ,n}:n\in \omega \}\cup \{x_{\gamma ,\ell }:\ell <d\}$ and $%
K_{\gamma }^{\prime }$ has a basis $\{w_{\gamma ,n}^{\prime }:n\in \omega \}$
where 
\begin{equation*}
w_{\gamma ,n}^{\prime }=p_{\gamma ,n}y_{\gamma ,n+1}-y_{\gamma
,n}-\sum_{\ell <d}s_{\gamma ,n,\ell }x_{\gamma ,\ell }
\end{equation*}
for some $p_{\gamma ,n}$, $s_{\gamma ,n,\ell }\in R$ where the $p_{\gamma
,n} $ are non-units of $R$ (not necessarily prime). (Compare, for example, 
\cite[Observation 3.1]{GS}.)

Let $F=\oplus _{\beta <\omega _{1}}F_{\beta }$ and $K=\oplus _{\beta <\omega
_{1}}K_{\beta }$ be as in \cite[Lemma XII.1.4]{EM}; that is, for all $\alpha
<\omega _{1}$, $\oplus _{\beta <\alpha }F_{\beta }/\oplus _{\beta <\alpha
}K_{\alpha }\cong A_{\alpha }$ and $\oplus _{\beta \leq \alpha }F_{\beta
}/(\oplus _{\beta <\alpha }F_{\beta }+K_{\alpha })\cong A_{\alpha
+1}/A_{\alpha }$. Moreover, by the proof of \cite[Lemma XII.1.4]{EM}, we can
assume that for $\gamma \in S$, $F_{\gamma }^{\prime }$ is a summand of $%
F_{\gamma }$ and $K_{\gamma }$ has a basis which includes $\{w_{\gamma
,n}:n\in \omega \}$ where 
\begin{equation*}
w_{\gamma ,n}=w_{\gamma ,n}^{\prime }-a_{\gamma ,n}
\end{equation*}
for some $a_{\gamma ,n}\in \oplus _{\beta <\gamma }F_{\beta }$ (and $\psi
_{\gamma }(w_{\gamma ,n}^{\prime })=\varphi _{\gamma }(a_{\gamma ,n})\in
A_{\gamma }$). Fix a basis $\mathcal{B}$ of $F$ which is the union of a
basis $\mathcal{B}_{\beta }$ for each $F_{\beta }$ and which includes $%
\bigcup_{\gamma \in S}\{y_{\gamma ,n}:n\in \omega \}\cup \{x_{\gamma ,\ell
}:\ell <d\}$. Also fix a basis of $K$ which includes $\bigcup_{\gamma \in
S}\{w_{\gamma ,n}:n\in \omega \}$. Given an element $r$ of $R$, we will say $%
\mu \in \omega _{2}$\textit{\ occurs in }$r$ if $r$ does not belong to $\mathbb{%
Q}(\{t_{\nu }:\nu \in \omega _{2}-\{\mu \}\})[X]$. Given an element $z$ of $%
F $ we will say that $\mu $\textit{\ occurs in }$z$ if it occurs in some
coefficient of the unique linear combination of elements of $\mathcal{B}$
which equals $z$. There is a subset $I$ of $\omega _{2}$ of cardinality $%
\aleph _{1}$ such that all of the $p_{\gamma ,n}$ and $s_{\gamma ,n,\ell }$ (%
$\gamma \in S$, $n\in \omega $, $\ell <d$) belong to $\mathbb{Q}(\{t_{i}:i\in
I\})[X]$. Moreover, we can choose $I$ such that it contains every $\mu $
which occurs in some coefficient of a linear combination of elements of $%
\mathcal{B}$ which equals some $a_{\gamma ,n}$ ($\gamma \in S$, $n\in \omega 
$). Without loss of generality (by renumbering the $t_{\nu }$), $I=\omega
_{1}$.

Now we define $\psi :K\rightarrow R$ by defining 
\begin{equation*}
\psi (w_{\gamma ,n})=t_{\omega _{1}+\omega \gamma +n}
\end{equation*}
and letting $\psi $ be arbitrary on the other basis elements of $K$. We will
show that $\limfunc{Ext}(A,R)\neq 0$ by showing that $\psi $ cannot be
extended to a homomorphism from $F$ into $R$. Suppose to the contrary that
there is a homomorphism $\varphi :F\rightarrow R$ extending $\psi $. For
each $\alpha <\omega _{1}$, let $T_{\alpha }$ be the set of all $\mu \in
\omega _{2}$ which occur in $\varphi (b)$ for some $b\in \bigcup \{\mathcal{B%
}_{\beta }:\beta <\alpha \}$. Then the $T_{\alpha }$ ($\alpha \in \omega _{1}
$) form a continuous chain of countable subsets of $\omega _{2}$ and there
is $\delta \in S$ such that $T_{\delta }\cap \{\omega _{1}+\beta :\beta
<\omega _{1}\}\subseteq \{\omega _{1}+\beta :\beta <\delta \}.$ There is a
finite subset $Z$ of $\omega _{2}$ such that every $\mu $ which occurs in $%
\varphi (y_{\delta ,0})$ or in $\varphi (x_{\delta ,\ell })$ for some $\ell
<d$ belongs to $Z$. Let $R^{*}=\mathbb{Q}(\{t_{\nu }:\nu \in \omega _{1}\cup
T_{\delta }\cup Z\})[X]$, a subring of $R=\mathcal{F}[X]$. Now for all $n\in
\omega $ we have $\varphi (w_{\delta ,n})=\psi (w_{\delta ,n})=$%
\begin{equation*}
t_{\omega _{1}+\omega \delta +n}=p_{\delta ,n}\varphi (y_{\delta
,n+1})-\varphi (y_{\delta ,n})-\sum_{\ell <d}s_{\delta ,n,\ell }\varphi
(x_{\delta ,\ell })-\varphi (a_{\delta ,n})\text{.}
\end{equation*}
If we can show that this implies that $t_{\omega _{1}+\omega \delta +n}$
belongs to $R^{*}$ for all $n\in \omega $, we will have a contradiction of
the choice of $T_{\delta }$ and the fact that $Z$ is finite. We will show
this by induction on $n$ along with simultaneously proving that $\varphi
(y_{\delta ,n+1})\in R^{*}$. We begin with $n=-1$: $\varphi (y_{\delta ,0})$
belongs to $R^{*}$ by definition of $Z$. Now suppose the inductive
hypothesis is true for $n-1$ and we prove it for $n$. By the last displayed
formula, the inductive hypothesis and the choice of $R^{*}$, there is an
element $r_{n}\in R^{*}$ such that $p_{\delta ,n}\varphi (y_{\delta
,n+1})=r_{n}-t_{\omega _{1}+\omega \delta +n}$. If $t_{\omega _{1}+\omega
\delta +n}\notin R^{*}$, there is an automorphism $\Theta $ of $R$ which
fixes $R^{*}$ and takes $t_{\omega _{1}+\omega \delta +n}$ to $t_{\tau }$
for some $\tau \notin T_{\delta }$. Then $p_{\delta ,n}\Theta (\varphi
(y_{\delta ,n+1}))=r_{n}-t_{\tau }$. (Remember that $p_{\delta ,n}\in R^{*}$%
.) Therefore, subtracting, $p_{\delta ,n}$ divides $t_{\omega _{1}+\omega
\delta +n}-t_{\tau },$which is impossible since $p_{\delta ,n}$ is a
non-unit. Thus $t_{\omega _{1}+\omega \delta +n}$ and hence $p_{\delta
,n}\varphi (y_{\delta ,n+1})$ belong to $R^{*}$. But then since $p_{\delta
,n}\in R^{*}$ we can prove by induction on $m$ that the coefficient of $X^{m}
$ in $\varphi (y_{\delta ,n+1})\in \mathcal{F}[X]$, belongs to $\mathbb{Q}%
(\{t_{\nu }:\nu \in \omega _{1}\cup T_{\delta }\cup Z\})$, and hence that $%
\varphi (y_{\delta ,n+1})$ belongs to $R^{*}$. \qed  

\medskip

We can even find a principal ideal domain of cardinality $\aleph _{1}$ which
satisfies the conclusion of Theorem \ref{main2}. Namely, let $R=\mathcal{F}%
_{1}[X]$ where $\mathcal{F}_{1}=\mathbb{Q}(\{t_{\nu }:\nu <\omega _{1}\})$.
Define $\psi (w_{\delta ,n})$ to be $t_{\omega \delta +\sigma _{\delta }+n}$
where $\omega \delta +\sigma _{\delta }$ is larger than any $\mu $ which
occurs in any $p_{\delta ,k}$ or $s_{\delta ,k,\ell }$ for $k\in \omega $, $%
\ell <d$. Define $T_{\alpha }$ as before and choose $\delta \in S$ such that 
$T_{\delta }\cap \omega _{1}\subseteq \omega \delta $. Let $R^{*}=\mathbb{Q}%
(\{t_{\nu }:\nu \in \omega \delta +\sigma _{\delta }\cup T_{\delta }\cup
Z\})[X]$.

We can also localize without affecting the property of the ring that we
desire. More generally, we have:

\begin{theorem}
\label{2b}For any $\kappa \geq \aleph _{1}$ there is a local p.i.d.\ $R$ of
cardinality $\kappa $ such that an $R$-module of cardinality $\leq \kappa $
is Whitehead only if it is free. \qed  
\end{theorem}

\end{document}